\documentclass{article}
\usepackage{latexsym,amsfonts,amssymb}
\usepackage[latin1]{inputenc}

\usepackage{tikz}

\usepackage{amsmath}

\title{\sc Graphs associated to conjugacy classes of normal subgroups in finite groups}

\author{Antonio Beltr\'an\\
\footnotesize
Departamento de Matem\'aticas,\\
\footnotesize Universidad Jaume I, \footnotesize
12071 Castell\'on, Spain\\
\footnotesize
e-mail: abeltran@mat.uji.es\\
\\
Mar\'{\i}a Jos\'e Felipe\\
\footnotesize
Instituto Universitario de Matem\'atica Pura y Aplicada,\\
\footnotesize Universidad Polit\'ecnica de Valencia, \footnotesize
46022 Valencia, Spain\\
\footnotesize e-mail: mfelipe@mat.upv.es \\
\\ Carmen Melchor\\
\footnotesize
Departamento de Educaci\'on,\\
\footnotesize Universidad Jaume I, \footnotesize
12071 Castell\'on, Spain\\
\footnotesize
e-mail: cmelchor@uji.es      }

\date{}
\begin{document} \maketitle

\begin{abstract}
Let $G$ be a finite group and let $N$ be a normal subgroup of $G$. We attach to $N$ two graphs $\Gamma_{G}(N)$ and $\Gamma_{G}^{*}(N)$ related to the conjugacy classes of $G$ contained in $N$ and to the set of primes dividing the sizes of these classes, respectively. These graphs are subgraphs of the ordinary ones associated to the conjugacy classes of $G$, $\Gamma (G)$ and $\Gamma^{*}(G)$, which have been widely studied by several authors. We prove that the number of connected components of both graphs is at most 2, we determine the best upper bounds for the diameters and characterize the structure of $N$ when these graphs are disconnected.



\end{abstract}

\bigskip

\section{Introduction}

Let $G$ be a finite group and let $N$ be a normal subgroup of $G$. For each element $x \in N$, the $G$-conjugacy class is $x^{G}=\lbrace x^{g}\, |\, g \in G \rbrace$. We will denote by Con$_{G}(N)$ the set of conjugacy classes in $G$ of elements of $N$. The elements in  Con$_{G}(N)$ are unions of conjugacy classes of $N$, and it turns out that every $G$-class size is a  multiple of a $N$-class size. Recent results have showed that the $G$-class sizes still have a strong influence on the structure of $N$ in spite of the fact that there may exist primes dividing the $G$-class sizes which, however, do not divide the order of $N$.

\bigskip

In 1990, E. A. Bertram, M. Herzog and A. Mann introduced in \cite{BerHerMann} the graph $\Gamma(G)$ associated to the sizes of the ordinary conjugacy classes of $G$, and later, in  \cite{Chi} the best bound of the diameter of this graph was attained.  Our aim is to study the properties of the following subgraph of $\Gamma(G)$ regarding the $G$-conjugacy classes contained in $N$ and to obtain structural properties of $N$ in the disconnected case.

\bigskip

{\bf Definition 1.1} Let $G$ be a finite group and let $N$ be a normal subgroup in $G$. We define the graph $\Gamma_{G}(N)$ in the following way: the set of vertices is the set of non-central elements of Con$_{G}(N)$, and two vertices $x^{G}$ and $y^{G}$ are joined by an edge if and only if $|x^{G}|$ and $|y^{G}|$ have a common prime divisor.
\bigskip

Notice that $\Gamma(N)$ is not a subgraph of $\Gamma_{G}(N)$ because the set of vertices of $\Gamma(N)$ needs not to be included within the set of vertices of $\Gamma_{G}(N)$.
 Moreover, we remark that although $\Gamma_{G}(N)$ is subgraph of $\Gamma(G)$, the fact that the number of connected components and the diameter of $\Gamma(G)$  are bounded does not directly imply that the corresponding for $\Gamma_{G}(N)$ have to be bounded too.  However, we show that both numbers, denoted by $n(\Gamma_{G}(N))$ and $d(\Gamma_{G}(N))$, are actually bounded. It is easy to check  that the bounds in Theorems A and B are the best possible bounds.

\bigskip

{\bf Theorem A}.
{\it Let $G$ a finite group and let $N$ be a normal subgroup of $G$. Then $n(\Gamma_{G}(N))\leq 2$.}

\bigskip
We want to remark that there is no relation between the connectivity of $\Gamma_G(N)$ and $\Gamma(N)$.  For instance, $\Gamma(N)$ can be disconnected while $\Gamma_G(N)$ is not. We can use the semilinear affine group $\Gamma(p^n)$ for appropriate $p$ and $n$ in order to see this. Recall that if $GF(p^{n})$ is the finite field of $p^{n}$ elements, then the multiplicative group $H=GF(p^{n})^{*}$ is cyclic of order $p^{n}-1$ and acts on the elementary abelian (additive) $p$-group of $GF(p^{n})$, say $K$. This action is Frobenius, so the corresponding semidirect product $KH$ is a Frobenius group with abelian kernel and complement. Moreover, $\alpha$, defined by $x^{\alpha}=x^{p}$ for all $x \in K$, is an automorphism of $K$ of order $n$ in such a way that $H\langle \alpha \rangle\leq $ Aut$(K)$. Then $\Gamma(p^n)$ is defined as the semidirect product $K(H\langle \alpha \rangle)$. Now, take $n=2$ and let $S$ be a cyclic subgroup of $H$ of order $s=3$ (so we are assuming that 3 divides $p^{2}-1$). We have that $N:=KS$ is normal in $G$ and is also a Frobenius group with abelian kernel and complement. Hence $\Gamma(N)$ is disconnected  by Theorem 2 of \cite{BerHerMann}. However, there are exactly two non-trivial $G$-classes in $N$ consisting in the $p^{2}-1$ elements of $K\setminus \lbrace 1 \rbrace$, and the $(|S|-1)|K|=2p^{2}$ elements of $NS \setminus K$, respectively.  Therefore, $\Gamma_G(N)$ is connected.

 \bigskip
 {\bf Theorem B}.
 {\it  Let $G$ a finite group and let $N$ be a normal subgroup of $G$.
\begin{enumerate}
\item  If $n(\Gamma_{G}(N))=1$, then ${d(\Gamma_{G}(N))\leq 3}$.
\item If $n(\Gamma_{G}(N))=2$, then each connected component is a complete graph.
\end{enumerate}
 }
\bigskip

We notice that the diameters of $\Gamma_G(N)$ and $\Gamma(N)$ are not either related. For instance, let $P$ be an extraspecial group of order $p^3$ with $p\neq 2$. If  we take $G=P\times S_3$ and $N=P\times A_3$, we have that $\Gamma(N)$ is a complete graph (all $N$-classes have size $p$) while $\Gamma_G(N)$ has diameter 2, since the nontrivial $G$-classes of $N$ have size $2, p$ and $2p$.

\bigskip
In 1995, S. Dolfi introduced in \cite{Dol} the dual graph $\Gamma^{*}(G)$  (it was also independently studied in \cite{Alf}) associated to the primes that divide the sizes of the conjugacy classes of $G$. In a similar way, we define the following subgraph of $\Gamma^{*}(G)$.

\bigskip

{\bf Definition 1.2} Let $G$ be a finite group and let $N$ be a normal subgroup in $G$. We define the ``dual" graph of $\Gamma_{G}(N)$, denoted by $\Gamma_{G}^{*}(N)$, as follows: the vertices are those primes which divide the size of some class in Con$_G(N)$,  and two vertices $p$ and $q$ are joined by an edge if there exists $C\in$ {\rm Con}$_{G}(N)$ such that $pq$ divides $\vert C\vert$.

\bigskip
 We also provide the best bounds for the number of components of $\Gamma^*_G(N)$ and for its diameter. We note again that these bounds cannot be obtained from the only fact that $\Gamma^*_G(N)$ is a subgraph of $\Gamma^*(G)$.

\bigskip
{\bf Theorem C}. {\it If $G$ is a finite group and $N\unlhd G$, then  $n(\Gamma_{G}^{*}(N))\leq 2$ and $n(\Gamma_{G}^{*}(N))=n(\Gamma_{G}(N))$.}

\bigskip

{\bf Theorem D}.  {\it  Let $G$ a finite group and $N\unlhd G$.

\begin{enumerate}
\item  If $n(\Gamma_{G}^{*}(N))=1$, then ${d(\Gamma_{G}^{*}(N))\leq 3}$.
\item  If $n(\Gamma_{G}^{*}(N))=2$, then each connected component is a complete graph.
\end{enumerate}
 }

 We give a characterization of  a normal subgroup $N$ whose graph $\Gamma_G(N)$, or equivalently $\Gamma_G^*(N)$, is disconnected. We recall that a group $G$ is said to be quasi-Frobenius if $G/{\bf Z}(G)$ is a Frobenius group. In this case, the inverse image in $G$ of the kernel and complement of $G/{\bf Z}(G)$ are called the kernel and complement of $G$, respectively.

\bigskip

{\bf Theorem E}.
{\it Let $G$ a finite group and $N\unlhd G$. If $\Gamma_{G}(N)$ has two connected components then, either $N$ is quasi-Frobenius with abelian kernel and complement, or $N=P\times A$ where $P$ is a $p$-group and $A\leqslant$ {\rm \textbf{Z}}$(G)$.}

\bigskip
We point out that our proofs of Theorems B, D and E are different and independent of the proofs of the respective theorems  concerning $\Gamma(G)$ and $\Gamma^*(G)$. In addition, Theorem E extends Corollary B of \cite{Akh}, which analyzes (with a completely different approach) the particular case in which $N$ has exactly two coprime $G$-class sizes bigger than $1$.\\

All groups considered will be finite and, if $A$ is a group or a set, $\pi (A)$ denotes the set of primes dividing $|A|$.

\section{Number of connected components of $\Gamma_G(N)$ and $\Gamma_G^*(N)$}

In this section we prove Theorems A and C in an easy way by using the following lemma, that is basic for our development. The distance in both graphs will be denoted by $d$.

\bigskip
{\bf Lemma 2.1} {\it Let $G$ be a finite group and $N\unlhd G$. Let $B=b^{G}$ y $C=c^{G}$ non-central elements in {\rm Con}$_{G}(N)$. If $(|B|,|C|)=1$. Then
\begin{enumerate}
\item {\rm \textbf{C}}$_{G}(b)${\rm \textbf{C}}$_{G}(c)=G$.
\item $BC=CB$ is a non-central element of {\rm Con}$_{G}(N)$ and $|BC|$ divides $|B| |C|$.
\item Suppose that $d(B,C)\geq3$ and $|B|<|C|$. Then $|BC|=|C|$ and $CBB^{-1}=C$. Furthermore, $C\langle BB^{-1}\rangle=C$, $\langle BB^{-1}\rangle\subseteq\langle CC^{-1}\rangle$ and $|\langle BB^{-1}\rangle |$ divides $|C|$.
\end{enumerate}}

\bigskip

{\it Proof.} For 1, 2 and the first part of 3, it is enough to mimic the proofs of Lemmas 1  and 2 of \cite{BerHerMann}, by taking into account that the product of two classes of Con$_{G}(N)$ is contained in $N$ again. The properties $C\langle BB^{-1}\rangle=C$ and $\langle BB^{-1}\rangle\subseteq\langle CC^{-1}\rangle$  are elementary. The fact that $|\langle BB^{-1}\rangle |$ divides $|C|$ follows from the fact that $C$ is a normal subset that can be written as the union of  right coclasses of the normal subgroup $\langle BB^{-1}\rangle$. $\Box$

\bigskip

{\it Proof of Theorem A.} Suppose that $\Gamma_{G}(N)$ has at least three connected components and take three non-central classes $B=b^{G}$, $C=c^{G}$ and $D=d^{G}$ in {\rm Con}$_{G}(N)$, each of which belongs to a different connected component. Certainly, any two of them have coprime size. Moreover, we can assume without loss of generality that $|B|< |C|< |D|$. By applying  Lemma 2.1, we get that $|\langle BB^{-1}\rangle|$ divides both $|D|$ and $|C|$. Then, $(|C|,|D|)>1$, which is a contradiction. $\Box$

\bigskip
{\it Proof of Theorem C.} Suppose that $n(\Gamma_{G}^{*}(N))\geq 3$. We take three primes $r, s$ and $l$ each of which belongs to a different connected component, and let $B, C$ and $D$ be  elements of Con$_{G}(N)$ such that $r$ divides $|B|$, $s$ divides $|C|$ and $l$ divides $|D|$. Without loss of generality we suppose that $|D|< |C| < |B|$. We have $d(B,D)\geq 3$ and $d(B,C)\geq 3$ and by applying Lemma 2.1, we obtain that $|\langle DD^{-1}\rangle|$ divides $|B|$ and $|C|$, but this leads to a contradiction, because $|B|$ and $|C|$ would have a common prime divisor. This proves that $n(\Gamma_{G}^{*}(N))\leq 2$. \\

Suppose now that $n(\Gamma_{G}(N))= 1$ and $n(\Gamma_{G}^{*}(N))=2$. Let $r$ and $s$ be primes such that each of them belongs to a distinct connected component of $\Gamma_{G}^{*}(N)$. Then there exist $B_{r}, B_{s} \in \Gamma_{G}(N)$ such that $r$ divides $|B_{r}|$ and $s$ divides $|B_{s}|$. Let us consider the following path in $\Gamma_{G}(N)$ that joins $B_{r}$ and $B_{s}$, which exists because $n(\Gamma_{G}(N))= 1$:
\begin{center}
\begin{tikzpicture}[node distance=1.5cm, auto]
  \node (A) {$B_{r}$};
  \node (B) [right of=A] {$B_{1}$};
  \node (C) [right of=B] {$B_{2}$};
  \node (D) [right of=C] {$\ldots$};
  \node (E) [right of=D] {$B_{s}$};
\draw[<->] (A) to node {$p_{1}$} (B);
\draw[<->] (B) to node {$p_{2}$} (C);
\draw[<->] (C) to node {$p_{3}$} (D);
\draw[<->] (D) to node {$p_{s}$} (E);
\end{tikzpicture}
\end{center}
where $p_i$ is a prime dividing $|B_i|$. This provides a contradiction, because $r$ and $s$ are connected in $\Gamma_{G}^{*}(N)$ by the following path:
\begin{center}
\begin{tikzpicture}[node distance=1.5cm, auto]
  \node (A) {$r$};
  \node (B) [right of=A] {$p_{1}$};
  \node (C) [right of=B] {$p_{2}$};
  \node (D) [right of=C] {$\ldots$};
  \node (E) [right of=D] {$p_{s}$};
   \node (F) [right of=E] {$s$};
\draw[<->] (A) to node {$B_{r}$} (B);
\draw[<->] (B) to node {$B_{1}$} (C);
\draw[<->] (C) to node {$B_{2}$} (D);
\draw[<->] (D) to node {$B_{s-1}$} (E);
\draw[<->] (E) to node {$B_{s}$} (F);
\end{tikzpicture}
\end{center}
So, we have proved that $n(\Gamma_{G}(N))= 1$ implies that $n(\Gamma_{G}^{*}(N))=1$. Now, if $n(\Gamma_{G}(N))= 2$ and $n(\Gamma_{G}^{*}(N))=1$ we can get a contradiction by arguing in a similar way. This shows that $n(\Gamma_{G}(N))=n(\Gamma_{G}^{*}(N))$. $\Box$

\section{Diameter of $\Gamma_G(N)$}

The following two lemmas, one for the disconnected case and the other for the connected case, summarize important structural properties of a normal subgroup $N$ concerning  the graph $\Gamma_G(N)$, which will be used for determining the diameters of $\Gamma_G(N)$ and $\Gamma^*_G(N)$. We start with the disconnected case.

 \bigskip

{\bf Lemma 3.1} {\it Let $G$ a finite group and let $N$ be a normal subgroup of $G$. Suppose that $n(\Gamma_{G}(N))=2$ and let $X_{1}$ and $X_{2}$ be the connected components of $\Gamma_{G}(N)$. Let $B_{0}$ be a non-central element of {\rm Con}$_{G}(N)$ of maximal size and assume that $B_{0} \in X_{2}$. We define
\begin{center}
$S=\langle C$ $|$ $C \in X_{1} \rangle$ and $T=\langle CC^{-1}$ $|$ $C \in X_{1} \rangle$.
\end{center}
Then
\begin{enumerate}

\item $S$ is a normal subgroup of $G$ and every element in S, either is central, or its $G$-conjugacy class is in $X_{1}$.
\item If $C$ is a $G$-conjugacy class of $N$ out of $S$, then $|T|$ divides $|C|$.
\item  $T=[S,G]$  is normal in $G$ and $T\leq $ {\rm\textbf{Z}}$(S)$.
\item {\rm \textbf{Z}}$(G)\cap N \subseteq S$ and $\pi(S/(${\rm \textbf{Z}}$(G)\cap N))\subseteq \pi(T)\subseteq \pi(B_{0})$. Moreover, $S$ is abelian.
\item Let $b^{G}=B \in X_{1}$. Then {\rm \textbf{C}}$_{G}(b)/S$ is a $q$-group for some prime $q \in \pi(B_{0})$.
\end{enumerate}

}

\bigskip

{\it Proof.} 1. The fact that $S$ is normal in $G$ is elementary. Let $C \in X_{2}$ and $B \in X_{1}$. We know that $BC$ is a $G$-conjugacy class of Con$_{G}(N)$ of maximal size between $|B|$ and $|C|$ by Lemma 2.1. Assume that $|BC|=|B|$. By Lemma 2.1 again, it follows that $|\langle CC^{-1}\rangle|$ divides $|B|$ and that $\langle CC^{-1}\rangle \subseteq \langle BB^{-1}\rangle$. On the other hand, $|B_0B|=|B_0|$ again by Lemma 2.1, and also $|\langle BB^{-1}\rangle|$ divides $|B_0|$. From these facts, we deduce that $(|B|, |B_{0}|)>1$, which is a contradiction. Thus, $|BC|=|C|$ for all $C\in X_{2}$ and $B \in X_{1}$. Furthermore, we have proved that the size of every class in $X_1$ is less than the size of any class in $X_2$.

 Now, take $C\in X_2$ and let $A$ be the union of all $G$-conjugacy classes of size $|C|$ in $S$ and assume that $A\neq \emptyset$. By the above paragraph, we have that if $B\in X_{1}$, then $BA\subseteq A$. Hence, $SA=A$,  and consequently, since $A$ is a normal subset,  $|S|$ divides $|A|$. This is not possible because $A \subseteq S-\lbrace 1\rbrace$. This contradiction shows that $A=\emptyset$, that is, $S$ does not contain any class of size $|C|$. Therefore, since $S$ is normal in $G$, then $S$ does not contain elements whose classes are in $X_2$.\\

2. Let $B\in X_1$. As we have proved in (1), $|B|<|C|$ and then, by Lemma 2.1, we have $C\langle BB^{-1}\rangle=C$, and as a consequence, $CT=C$. Therefore, $|T|$ divides $|C|$, as wanted.\\

3. By definition, it is clear that $T=[S,G]$ and so, it is a normal subgroup of $G$. Let us prove that $T\leq $ {\rm \textbf{Z}}$(S)$. In fact, if $B=b^{G} \in X_{1}$, then $(|T|, |G:{\rm \textbf{C}}_{G}(b)|)=(|T|, |B|)=1$, because $|T|$ divides every class size in $X_{2}$ by (2). Now, since $|T : {\rm \textbf{C}}_{T}(b)|$ divides $(|T|, |B|)=1$, we deduce that $T = {\rm \textbf{C}}_{T}(b)$. As the classes in $X_{1}$ generate $S$, we conclude that $T$ is central in $S$. \\

4.  Let $z \in \textbf{Z}(G)\cap N $ and let $B=b^{G}\in X_1$. Note that $b^{G}z=(bz)^{G}$. Moreover $bz \in N$, because both elements lie in $N$. As $|(bz)^{G}| =|Bz|=|B|$, then $bz \in S$ and so $z\in S$. This proves that $\textbf{Z}(G)\cap N \subseteq S$.   Since $T=[S,G]$, then $[S/T,G/T]=1$ and $S/T\subseteq$ \textbf{Z}$(G/T)$. In particular, $S/T$ is abelian and as $T\leq {\bf Z}(S)$ by (3), then $S$ is nilpotent. We can write $S=R \times Z$ where $Z$ is the largest Hall subgroup of $S$ which is contained in \textbf{Z}($G$). Let $p$ be a prime divisor of $|R|$ and let $P$ be a Sylow $p$-subgroup of $R$. It is clear that $P \unlhd G$ and $T=[S, G]=[R, G]\geqslant [P, G]> 1$. Hence $p$ divides $|T|$ and by applying (1) and (2),  $|T|$ divides $|B_0|$. Therefore, $\pi(R)\subseteq \pi(T) \subseteq \pi(B_{0})$. On the other hand,   it is elementary that $\pi(S/(${\rm \textbf{Z}}$(G)\cap N))\subseteq \pi(R)$, and the first part of the step is proved. We show now that  $R\leq {\rm \textbf{Z}}(S)$. In fact, let $b^G= B \in X_{1}$. Since $(|B|, |B_{0}|)=1$, we obtain in particular, $(|B|,|R|)=1$. Thus, $|R:$\textbf{C}$_{R}(b)|=1$ since this index trivially divides $|R|$ and $|B|$ because $R\unlhd G$. This means that $R=$ \textbf{C}$_{R}(b)$ for every generating element $b$ of $S$. So, $R$ is contained in \textbf{Z}$(S)$ as wanted, and $S$ is abelian. \\

5.  By considering the primary decomposition of $b$, it is clear that we can write $b=b_{q}b_{q'}$ where $b_{q}$ and $b_{q'}$ are the $q$-part and the $q'$-part of $b$, where $q$ is a prime such that $b_{q}\not\in \textbf{Z}(G)\cap N$. Hence, $q \in  \pi(B_{0})$ by (4). Furthermore, it is elementary that $\textbf{C}_{G}(b)\subseteq \textbf{C}_{G}(b_{q})$, and as a result, $|(b_{q})^{G}|$ divides $|B|$. We claim that any element $xS \in \textbf{C}_{G}(b)/S$ is a $q$-element. For any $x\in  \textbf{C}_{G}(b)$, write $x=x_{q}x_{q'}$ (it is possible $x_{q}=1$). It is obvious that $x_{q}$ and $x_{q'}$ belong to \textbf{C}$_{G}(b)$. We consider $a=b_{q}x_{q'}$ and observe that $\textbf{C}_{G}(a)=\textbf{C}_{G}(b_{q})\cap \textbf{C}_{G}(x_{q'})\subseteq \textbf{C}_{G}(b_{q})$, so $|(b_{q})^{G}|$ divides $|a^{G}|$. Since $(b_q)^G\in X_1$, this forces that $a^G\in X_1$, and  we conclude that $x_{q'}\in S$, that is, $xS$ is a $q$-element, as wanted. This shows that {\rm \textbf{C}}$_{G}(b)/S$ is a $q$-group. $\Box$

\bigskip

{\bf Lemma 3.2} {\it Let $G$ be a finite group and $N\unlhd G$ with $\Gamma_{G}(N)$ connected. Let $B_{0}$ be a $G$-conjugacy class of {\rm Con}$_{G}(N)$ of maximal size. Let
$$M=\langle D \mid D \in {\rm Con}_{G}(N)\, \, and \,\,  d(B_{0},D)\geq 2\rangle$$
$$K=\langle D^{-1}D \mid D \in {\rm Con}_{G}(N)\,\, and\, \, d(B_{0},D)\geq 2\rangle$$
Then
\begin{enumerate}
\item $M$ and $K$ are normal subgroups of $G$. Furthermore, $K=[M,G]$ and $K\leq $ {\rm \textbf{Z}}$(M)$.
\item {\rm \textbf{Z}}$(G)\cap N \subseteq M$ and $\pi(M/(${\rm \textbf{Z}}$(G)\cap N))\subseteq \pi(K)\subseteq \pi(B_{0})$. Furthermore, $M$ is abelian.

\end{enumerate}
}

\bigskip
{\it Proof.}  1. By definition, we easily see that $M$ and $K$ are normal subgroups of $G$ and $K=[M,G]$. Let us prove that $K\leq$ {\rm \textbf{Z}}$(M)$. If $C=c^G \in$ {\rm Con}$_{G}(N)$ satisfies that  $d(B_{0},C)\geq 2$, in particular we have $(|B_{0}|,|C|)=1$ and then, $|B_{0}|=|B_{0}C|$. Moreover, by Lemma 2.1,  $B_{0}CC^{-1}=B_{0}$ and as a result $|K|$ divides $|B_{0}|$. Therefore, $(|K|,|C|)=1$. However, we have that $|K:\textbf{C}_{K}(c)|$ divides $(|K|,|C|)$ and thus, $K=\textbf{C}_{K}(c)$, which implies that $K\leq$ {\rm \textbf{Z}}$(M)$.\\

2. We prove that $\textbf{Z}(G)\cap N \subseteq M$. Let $z \in \textbf{Z}(G)\cap N $ and let $C=c^{G}\in$ {\rm Con}$_{G}(N)$ such that $d(B_{0}, C)=2$. Notice that $c^{G}z=(cz)^{G}$. As $|(cz)^{G}| =|c^{G}|$, then $d(B_{0}, (cz)^{G})=2$. Thus, $cz \in M$ and $\textbf{Z}(G)\cap N \subseteq M$.  Since $K=[M,G]$, then $M/K\leq $ \textbf{Z}$(G/K)$ and since $K\leq $ \textbf{Z}$(M)$ by (1), we obtain that $M$ is nilpotent. We can write $M=R \times Z$ where $Z$ is the largest Hall subgroup of $M$ that is contained in $\textbf{Z}(G)$. Let $q$ be a prime divisor of $|R|$ and let $Q$ be the Sylow $q$-subgroup of $R$. Then $Q\unlhd G$ and $K=[M,G]\geqslant [R,G] \geqslant [Q, G]>1$. So, $q$ divides $|K|$ and $\pi(R)\subseteq \pi (K)$. In the proof of (1), we have seen that $\pi(K)\subseteq \pi(B_{0})$. Then, $\pi(R)\subseteq \pi( B_{0})$. Furthermore, it is elementary that $\pi(M/(${\rm \textbf{Z}}$(G)\cap N))\subseteq \pi(K)$ and so, the first part of the step is proved. On the other hand, given a generating class  $B=b^{G}$ of  $M$, we know that $d(B,B_{0})\geq 2$. In particular, we have $(|B|,|B_{0}|)=1$ and hence $(|Q|,|B|)=1$, where $Q$ is the above Sylow $q$-subgroup. Since $|Q:\textbf{C}_{Q}(b)|$ divides $(|Q|,|B|)=1$, we have $\textbf{C}_{G}(b)=Q$ and $Q\leq $ {\rm \textbf{Z}}$(M)$. Thus, $R\leq$ {\rm \textbf{Z}}$(M)$ and $M$ is abelian. $\Box$\\

The following consequence, which has interest on its own,  is the key to bound the diameter of $\Gamma_{G}(N)$ in the connected case.

\bigskip
{\bf Theorem 3.3}.
{\it Let $G$ a finite group and $N$ a normal subgroup of $G$ and suppose that  $\Gamma_{G}(N)$ is connected. Let $B_{0}$ a non-central conjugacy class of {\rm Con}$_{G}(N)$ with maximal size. Then $d(B,B_{0})\leq 2$ for every non-central $B \in$ {\rm Con}$_{G}(N)$. }

\bigskip
{\it Proof.} Suppose that the theorem is false. Let $B=b^{G}\in$ {\rm Con}$_{G}(N)$ such that $d(B_{0},B)=3$ and let
\begin{center}
\begin{tikzpicture}[node distance=1.5cm, auto]
  \node (A) {$B_{0}$};
  \node (B) [right of=A] {$B_{1}$};
  \node (C) [right of=B] {$B_{2}$};
  \node (D) [right of=C] {$B$};
\draw[<->] (A) to node {} (B);
\draw[<->] (B) to node {} (C);
\draw[<->] (C) to node {} (D);
\end{tikzpicture}
\end{center}
be a shortest chain linking $B$ and $B_{0}$ of length 3.  By considering the primary decomposition of $b$,  we  write $b=b_{q}b_{q'}$ where $b_{q}$ and $b_{q'}$ are the $q$-part and the $q'$-part of $b$, and $q$ is a prime such that $b_{q}\not\in \textbf{Z}(G)\cap N$. Hence, $q \in  \pi(B_{0})$ by Lemma 3.2(2). Also, the fact that $\textbf{C}_{G}(b)\subseteq \textbf{C}_{G}(b_{q})$ implies that $|(b_{q})^{G}|$ divides $|B|$ and then
 any class which is connected to $(b_{q})^{G}$ must be connected to $B$. This means that $d((b_{q})^{G},B_{0}))\geq 3$.

  Let $M$ be  subgroup defined in Lemma 3.2.  We claim that any element $x \in \textbf{C}_{G}(b_q) \setminus M$ satisfies that $xM$ is a $q$-element. Write $x=x_{q}x_{q'}$ and suppose that
       $x_{q'}\notin M$. Set $a=x_{q'}b_{q}$ and notice that $a \notin M$. By definition of $M$, we have $d(a^{G},B_{0})\leq 1$ and since $\textbf{C}_{G}(a)=\textbf{C}_{G}(x_{q'})\cap \textbf{C}_{G}(b_q)$, it follows that $|(b_q)^G|$ divides $|a^{G}|$. These facts show that $d((b_q)^G,B_{0}))\leq 2$, a contradiction. Therefore, $x_{q'}\in M$ and $xM$ is a $q$-element. In conclusion, $\textbf{C}_{G}(b_q)/M$ is a $q$-group.
  Now, observe that  $|B_{2}|$ divides
$$|G:\textbf{Z}(G)\cap N|=|G:\textbf{C}_{G}(b_q)|\, |\textbf{C}_{G}(b_q):M|\, |M:\textbf{Z}(G)\cap N|.$$
Also, we know by Lemma 3.2(2) that $\pi(M/\textbf{Z}(G)\cap N)\subseteq \pi(B_{0})$ and we have seen in the above paragraph that \textbf{C}$_{G}(b_q)/M$ is a $q$-group for some $q \in \pi(B_{0})$. Consequently, $|B_{2}|$ must divide $|(b_q)^G|$ (which divides $|B|$), because $(|B_2|,|B_0|)=1$. This is a contradiction, since $B_{1}$ and $B$ would be joined by an edge. $\Box$

\bigskip

{\it Proof of Theorem B.} 1. Suppose that $D_{1}$ and $D_{2}$ are classes of {\rm Con}$_{G}(N)$ such that $d(D_{1},D_{2})=4$. Let $B_{0}$ be a class of maximal size in {\rm Con}$_{G}(N)$. By Theorem 3.3 we know that $d(B_{0},D_{i})\leq 2$ for $i=1,2$. We can suppose then that $d(B_{0},D_{i})=2$ for $i=1,2$. Furthermore, without loss of generality, $|D_{1}|>|D_{2}|$. Then, by Lemma 2.1 it is true that $|\langle D_{2}D_{2}^{-1}\rangle |$ divides $|D_{1}|$. In addition, $B_{0}D_{2}$ is a conjugacy class of {\rm Con}$_{G}(N)$ and $|B_{0}D_{2}|=|B_{0}|$ by Lemma 2.1(2) and by the maximality of $B_0$. It follows that $B_0D_2D_2^{-1}=B_0$ and  $|\langle D_{2}D_{2}^{-1}\rangle |$ divides $|B_{0}|$. Therefore, $B_{0}$ and $D_{1}$ are joined by an edge, which is a contradiction. This proves that $d(\Gamma_G(N))\leq 3$.\\


2. Let $B_{1}=b_{1}^{G}$ and $B_{2}=b_{2}^{G}$ in $X_{1}$. Notice that $b_{1}, b_{2} \in S$ where $S$ is the subgroup defined in Lemma 3.1. By applying the properties of that result, we know that $|b_{2}^{G}|$ divides
$$|G:{\rm\textbf{Z}}(G)\cap N|=|G:{\rm \textbf{C}}_{G}(b_{1})||{\rm \textbf{C}}_{G}(b_{1}):S||S:{\rm\textbf{Z}}(G)\cap N|$$
where the primes dividing $|{\rm \textbf{C}}_{G}(b_{1}):S|$ and $|S:{\rm\textbf{Z}}(G)\cap N|$ are in $\pi(B_{0})$. So, we have that $|b_{2}^{G}|$ divides $|b_{1}^{G}|$. By arguing symmetrically we also get that $|b_{1}^{G}|$ divides $|b_{2}^{G}|$, so we conclude that all classes in $X_{1}$ have the same size. Hence, $X_{1}$ is a complete graph. Now, we prove that $X_{2}$ is also a complete graph.  It is enough to consider again $S$ and $T$  defined in Lemma 3.1 and observe that every $C \in X_{2}$ is out of $S$ and that $|T|$ divides $|C|$ by Lemma 3.1(1) and (2).
$\Box$




\bigskip
{\bf Remark 3.4}. In the proof of Theorem B.2, we have seen that all $G$-classes of $N$ lying in the connected component $X_1$ (the component which does not contains the classes of maximal size), must have the same size. Moreover, in the proof of Lemma 3.1(1) we have seen that this size is less than the size of every class in $X_2$.

\section{Diameter of $\Gamma_G^*(N)$}
{\it Proof of Theorem D.} 1. Suppose that there exist two primes $r$ and $s$ in $\Gamma_{G}^{*}(N)$ such that $d(r,s)=4$ and we will get a contradiction. This means that the primes $r$ and $s$ are connected by a path of length 4, say
\begin{center}
\begin{tikzpicture}[node distance=1.5cm, auto]
  \node (A) {$r$};
  \node (B) [right of=A] {$p_{1}$};
  \node (C) [right of=B] {$p_{2}$};
  \node (D) [right of=C] {$p_{3}$};
  \node (E) [right of=D] {$s$};
\draw[<->] (A) to node {$B_{1}$} (B);
\draw[<->] (B) to node {$B_{2}$} (C);
\draw[<->] (C) to node {$B_{3}$} (D);
\draw[<->] (D) to node {$B_{4}$} (E);
\end{tikzpicture}
\end{center}
where $B_{i}\in$ {\rm Con}$_{G}(N)$ for $i=1,\ldots,4$ and $p_{i}\in \Gamma_{G}^{*}(N)$ for $i=1, 2, 3$. By Theorem 3.3 we know that $d(B_{i},B_{0})\leqslant 2$ for $i=1,\ldots , 4$ where $B_{0}$ is a non-central $G$-conjugacy class of maximal size. Notice that $d(B_{1},B_{4})=3$ and we distinguish  only two possibilities:\\

\textit{Case 1.} $d(B_{0},B_{1})=2=d(B_{0},B_{4})$. By symmetry, we can assume for instance that $|B_{1}|>|B_{4}|$. Since $d(B_{1},B_{4})=3$, by Lemma 2.1 we have that $|\langle B_{4}B_{4}^{-1}\rangle |$ divides $|B_{1}|$. Moreover, $B_{0}B_{4}$ is an element of {\rm Con}$_{G}(N)$ such that $|B_{0}B_{4}|=|B_{0}|$ and by Lemma 2.1, $|\langle B_{4}B_{4}^{-1}\rangle |$ divides $|B_{0}|$. Therefore, $d(B_{0},B_{1})=1$, because their cardinalities have a prime common divisor. This is a contradiction.\\

\textit{Case 2.} Either $d(B_{0},B_{1})=2$ and $d(B_{0},B_{4})=1$, or $d(B_{0},B_{1})=1$ and $d(B_{0},B_{4})=2$. Without loss we assume for instance the latter case.  Let us consider the subgroup $M$ defined in Lemma 3.2 and  let $B_{4}=b^{G}$. Since $d(B_{0},B_{4})=2$, then $b \in M$ by definition. Moreover, $|B_{1}|$ divides
$$|G:\textbf{Z}(G)\cap N|=|G:\textbf{C}_{G}(b)| |\textbf{C}_{G}(b):M| |M:\textbf{Z}(G)\cap N|.$$
 Now, notice that $r \not \in \pi(B_{0})$, otherwise  it yields $d(r,s)\leq 3$, a contradiction,  and trivially $r\notin \pi(B_{4})$. Also, $\pi(M/\textbf{Z}(G)\cap N )\subseteq \pi(B_{0})$, so we have that $r$ (which divides $|B_1|$)  must divide $|{\rm \textbf{C}}_{G}(b):M|$. Therefore, there exists an $r$-element $y \in {\rm \textbf{C}}_{G}(b)\setminus M$. On the other hand,  $b\in M$, and by Lemma 3.2(2), the $r$-part of $b$ is central in $G$, that is,  we can assume that $b$ is an $r'$-element, by replacing $b$ by its $r$'-part. As $y$ and $b$ have coprime orders, we have
 $${\rm \textbf{C}}_{G}(yb)={\rm \textbf{C}}_{G}(y)\cap {\rm \textbf{C}}_{G}(b)\subseteq {\rm \textbf{C}}_{G}(b). $$
Consequently, $|B_4|$ divides $|(yb)^{G}|$. Furthermore, since $yb \not \in M$,  by definition of $M$ we have $d((yb)^{G},B_{0})\leq 1$. As $d(B_0,B_1)=1$ by hypothesis,  we deduce that $d(B_1,(yb)^G)\leq 2$. Now, $s$ divides $|B_1|$ and $r$ divides $|(yb)^G|$, and this forces that  $d(r,s)\leq 3$, which is a contradiction.\\

2. Let $X_{1}$ and $X_{2}$ be the connected components of $\Gamma_{G}(N)$ where $X_{2}$ is the component that contains the $G$-conjugacy class with the largest size. Let us prove first that $X_{1}^{*}$, $X_{2}^{*}$ are the connected components of $\Gamma_{G}^{*}(N)$, where $X_{i}^{*}=\lbrace p \in \pi(B)\,  |\, B \in X_{i} \rbrace$, and secondly, that $X_{1}^{*}$ y $X_{2}^{*}$ are complete graphs.\\

Let $X$ be a connected component of $\Gamma_{G}(N)$ and let $r,s \in X^{*}$. Then there exist $B_{r}, B_{s}$ such that $r$ divides $|B_{r}|$ and $s$ divides $|B_{s}|$. Let us consider de path in $\Gamma_{G}(N)$ that joins $B_{r}$ and $B_{s}$:
\begin{center}
\begin{tikzpicture}[node distance=1.5cm, auto]
  \node (A) {$B_{r}$};
  \node (B) [right of=A] {$B_{1}$};
  \node (C) [right of=B] {$B_{2}$};
  \node (D) [right of=C] {$\ldots$};
  \node (E) [right of=D] {$B_{s}$};
\draw[<->] (A) to node {$p_{1}$} (B);
\draw[<->] (B) to node {$p_{2}$} (C);
\draw[<->] (C) to node {$p_{3}$} (D);
\draw[<->] (D) to node {$p_{s}$} (E);
\end{tikzpicture}
\end{center}
where $B_{i}\in$ {\rm Con}$_{G}(N)$ for $i=1,\ldots,s-1$ and $p_{i}\in \Gamma_{G}^{*}(N)$ for $i=1,\ldots, s$. So, $r$ and $s$ are connected in $\Gamma_{G}^{*}(N)$ in the following way:
\begin{center}
\begin{tikzpicture}[node distance=1.5cm, auto]
  \node (A) {$r$};
  \node (B) [right of=A] {$p_{1}$};
  \node (C) [right of=B] {$p_{2}$};
  \node (D) [right of=C] {$\ldots$};
  \node (E) [right of=D] {$p_{s}$};
   \node (F) [right of=E] {$s$};
\draw[<->] (A) to node {$B_{r}$} (B);
\draw[<->] (B) to node {$B_{1}$} (C);
\draw[<->] (C) to node {$B_{2}$} (D);
\draw[<->] (D) to node {$B_{s-1}$} (E);
\draw[<->] (E) to node {$B_{s}$} (F);
\end{tikzpicture}
\end{center}
Therefore, $X^{*}$ is contained in a connected component $Y$ of  $\Gamma_{G}^{*}(N)$. Now, we take $q \in Y$, which is connected by an edge to some $r \in X^{*}$. Then there exists $B \in$ Con$_{G}(N)$ such that $qr$ divides $|B|$. It follows that $B \in X$ and $q \in X^{*}$. Thus, $X^{*}=Y$ and $X^{*}$ is a connected component of $\Gamma_{G}^{*}(N)$ as wanted.\\

 By Remark 3.4, all classes in $X_1$ have the same size, which trivially implies that $X_1^*$ is a complete graph. Let us show that $X_{2}^{*}$ is a complete graph too. Suppose that $B_{0}$ is a conjugacy class with maximal size, which lies in $X_{2}$, and let $B_{1}=b_{1}^{G}\in X_{1}$. Then, the subgroup $S$ defined in Lemma 3.1 is abelian, and $S \subseteq${\bf  C}$_{G}(b_{1})$. Now, if $p\in X_{2}^{*}$,  there exists $D \in X_{2}$ such that $p$ divides $|D|$. Notice that $|D|$ divides
$$|G:{\rm \textbf{Z}}(G)\cap N|=|B_{1}||\textbf{C}_{G}(b_{1}):S||S:\textbf{Z}(G)\cap N|,$$
and by Lemma 3.1(4) and (5), we know that $|$\textbf{C}$_{G}(b_{1}):S|$ is a $q$-power with $q\in \pi(B_{0})$ and $\pi(S/($\textbf{Z}$(G)\cap N))\subseteq \pi(B_{0})$. It follows that $\pi(D)\subseteq \pi(B_{0})$. Therefore, all primes in  $X_{2}^{*}$ are in  $\pi(B_{0})$ and so, $X_2^*$ trivially is a complete graph. $\Box$

\section{Structure of $N$ in the disconnected case}

{\it Proof of Theorem E.} Suppose that $X_1$ and $X_2$ are the two connected components of $\Gamma_G(N)$, where $X_2$ is the one containing the classes of maximal size. Let $S$ the subgroup defined in Lemma 3.1.\\

\textit{Step 1:} If $S\leq$ \textbf{Z}$(N)$, then $N=P\times A$ with $A\leq$ \textbf{Z}$(G) \cap N$ and $P$ a $p$-group.\\

We can choose  a $p$-element $x$ and a $q$-element $y$ of $N$, for some primes $p$ and $q$, such that $x^G\in X_{1}$ and $y^G \in X_2$. If $p=q$ for every election of $x$ and $y$, it is clear that $N = P \times A$ with $A \leq \textbf{\mbox{Z}}(G) \cap N$. Assume then that  $p\neq q$. Since $x \in S \leq \textbf{\mbox{Z}}(N)$, we obtain  $N/S = \textbf{C}_{N}(x)/S$ and this group has prime power order by Lemma 3.1(5). As a consequence, $N$ is nilpotent and $[x,y]=1$. As $x$ and $y$ have coprime order, $\textbf{C}_{G}(xy) = \textbf{C}_{G}(x) \cap \textbf{C}_{G}(y)$ and so, $|(xy)^G|$ divides $|x^G|$ and $|y^G|$, a contradiction. This proves the step. \\

Notice that we can assume that $S$\textbf{Z}$(N)< N$, because if $S$\textbf{Z}$(N)= N$ then $N$ is abelian and $S\leq$ \textbf{Z}$(N)$. For the rest of the proof, we assume \textbf{Z}$(N)<S$\textbf{Z}$(N)< N$ and we will prove that $N$ is quasi-Frobenius with abelian kernel and complement. We divide the proof into several steps. Let us denote by $\pi=\lbrace p$ prime $\vert$ $p$ divides $|B|$ with $B \in X_{1}\rbrace$.\\

\textit{Step 2:} $N$ has a normal $\pi$-complement and abelian Hall $\pi$-subgroups. \\

 Let us prove that $N$ is $p$-nilpotent and has abelian Sylow $p$-subgroups for every $p\in \pi$. Let $a \in N\setminus S$, then obviously $a^{G}\in X_{2}$ and $|a^{G}|$ is a $\pi'$-number. If we take $P \in $ Syl$_{p}(N)$, then there exists $g\in N$ such that $P^{g}\subseteq \textbf{C}_{N}(a)$, that is, $a \in$ \textbf{C}$_{N}(P^{g})=$\textbf{C}$_{N}(P)^{g}$. Thus, we can write $$N=S\cup \bigcup_{g\in N}{\rm\textbf{C}}_{N}(P)^{g}$$ and, by counting elements,  it follows that $$|N|\leq (|S|-1)+|N:{\rm \textbf{N}}_{N}({\rm \textbf{C}}_{N}(P))|(|{\rm \textbf{C}}_{N}(P)|-1)+1.$$
  Hence
  $$1\leq \frac{|S|}{|N|}+\frac{|{\rm\textbf{C}}_{N}(P)|}{|{\rm \textbf{N}}_{N}({\rm\textbf{C}}_{N}(P))|}- \frac{1}{|{\rm \textbf{N}}_{N}({\rm\textbf{C}}_{N}(P))|}.$$
However, if  \textbf{C}$_{N}(P)< \textbf{N}_{N}($\textbf{C}$_{N}(P))$, as we are assuming that $S<N$, we have $$1\leq \frac{1}{2}+\frac{1}{2}- \frac{1}{|{\rm \textbf{N}}_{N}({\rm\textbf{C}}_{N}(P))|},$$ which is a contradiction. This implies that \textbf{C}$_{N}(P)=\textbf{N}_{N}($\textbf{C}$_{N}(P))$, and in particular,  $$P\leq {\rm \textbf{N}}_{N}(P)\leq {\rm \textbf{N}}_{N}({\rm \textbf{C}}_{N}(P))\leq {\rm \textbf{C}}_{N}(P),$$
    so \textbf{C}$_{N}(P)=  \textbf{N}_{N}(P)$ and $P$ is abelian. By Burnside's $p$-nilpotency criterion (see for instance 17.9 of \cite{Hu}), we get that $N$ is $p$-nilpotent for every $p\in \pi$ and so, $N$ has normal $\pi$-complement. In particular, $N$ is $\pi$-separable and there exists a Hall $\pi$-subgroup $H$ of $N$. By reasoning  with $H$ similarly as with $P$, we obtain \textbf{C}$_{N}(H)=\textbf{N}_{N}(H)$ and so, $H$ is abelian too. The step is finished.\\

Let $K/$\textbf{Z}$(N)$ be the normal $\pi$-complement of $N/$\textbf{Z}$(N)$. By applying Lemma 3.1(4), we get that $S$\textbf{Z}$(N)/$\textbf{Z}$(N)$ is a normal $\pi'$-subgroup of $N$/\textbf{Z}$(N)$, so $S\leq K$.\\

\textit{Step 3:} $K = \textbf{C}_{N}(x)$ for every $x \in S\setminus$ \textbf{Z}$(G)\cap N$ and $S\leq$ \textbf{Z}$(K)$.\\

 Let $x\in S\setminus$ \textbf{Z}$(G)\cap N$. Then $x^{G}\in X_{1}$ by Lemma 3.1(1) and \textbf{C}$_{G}(x)/S$ is a $\pi'$-group by Lemma 3.1(5). Since $|x^N|$ is a $\pi$-number, we obtain that \textbf{C}$_{N}(x)/$\textbf{Z}$(N)$ is a Hall $\pi'$-subgroup of $N/$\textbf{Z}$(N)$. Thus, $K=$ \textbf{C}$_{N}(x)$ for every $x \in S\setminus$ \textbf{Z}$(G)\cap N$ and, in particular, $S\leq$ \textbf{Z}$(K)$.\\

\textit{Step 4:} $K = S$. \\

Let $H$ be an abelian Hall $\pi$-subgroup of $N$. We have seen in the proof of Step 2 that $$N=S\cup \bigcup_{g\in N}{\rm\textbf{C}}_{N}(H)^{g},$$ which trivially implies that $$N=\bigcup_{g\in N}S{\rm\textbf{C}}_{N}(H)^{g}.$$
This forces that $N=\textbf{C}_{N}(H)S$ and consequently,  $HS\unlhd N$.  Suppose that $S <K$ and we will get a contradiction.  Let $a \in K\setminus S$, then $a^{G}\in X_{2}$ by Lemma 3.1(1), so $|a^G|$ is a $\pi'$-number and as a result, $a \in$ \textbf{C}$_{K}(H^{g})=$ \textbf{C}$_{K}(H)^{g}$ for some $g \in N$. Moreover, $S\leq$ \textbf{Z}$(K)$ by Step 3, so we have the following equalities
$$\textbf{C}_{K}(H^{g})= \textbf{C}_{K}(H^{g}S)= \textbf{C}_{K}(HS)= \textbf{C}_{K}(H).$$
 Thus, $a \in \textbf{C}_{K}(H)$ for every $a \in K \setminus S$ and we conclude that $K = \langle K \setminus S \rangle \subseteq \textbf{C}_{N}(H)$. As $H$ is abelian and $N=HK$, we have $H\leq$ \textbf{Z}$(N) \leq K$. This implies that $N=K$ and then, $S \leq \textbf{Z}(N)$ by Step 3, which contradicts the assumption made after Step 1. \\

\textit{Step 5:} $N$ is quasi-Frobenius with abelian kernel and complement.\\

  Let $\overline{N}= N/\textbf{Z}(N)$ and let $\overline{K}= K/\textbf{Z}(N)$. If $\overline{K}=$ \textbf{C}$_{\overline{N}}(\overline{x})$ for all $\overline{x} \in \overline{K}\setminus \lbrace 1\rbrace$, this is equivalent to the fact that $N$ is quasi-Frobenius with abelian kernel $K$ and abelian complement $H$. Suppose then that $\overline{K}<$ \textbf{C}$_{\overline{N}}(\overline{x})$ for some $\overline{x} \in \overline{K}\setminus \lbrace 1\rbrace$. Also, we can suppose that $o(\overline{x})$ is an $r$-number for some prime $r \in \pi'$. Now, let $\overline{y}\in$C$_{\overline{N}}(\overline{x})\setminus \overline{K}$ such that $o(\overline{y})$ is a $q$-number for some $q \in \pi$. We can suppose without loss of generality that $o(y)$ is a $q$-number. Notice that $[y,x]\in$ \textbf{Z}$(N)$ because $\overline{y} \in$ \textbf{C}$_{\overline{N}}(\overline{x})$ and, since $(o(x),o(y))=1$, it easily follows that $[x,y]=1$. Then $y \in$ \textbf{C}$_{N}(x)=K$ and $\overline{y}\in \overline{K}$, which is a contradiction. $\Box$

\bigskip

{\bf Example 5.1}  We show that the converse of the above theorem is false. It is known that the special linear group $H=SL(2,5)$ acts Frobeniusly on $K\cong {\mathbb Z}_{11} \times {\mathbb Z}_{11}$. As a consequence, the action  of any subgroup of $H$ on $K$ is also Frobenius. We consider, in particular, a Sylow $5$-subgroup $P$ of $ H $ and $ {\rm \textbf{N}}_H(P)$ acting Frobeniusly on $K$. We define the semidirect product $N:= KP$, which is trivially a normal subgroup of $G:=K${\rm \textbf{N}}$_H(P)$. Thus, $N$ is a Frobenius group with abelian kernel and  complement. In fact, $N$ decomposes into the following disjoint union
$$N= \{1\} \cup  (K\setminus \{1\})  \cup  (\bigcup_{k\in K}  P^k\setminus \lbrace 1 \rbrace ), $$
and  $K\setminus \{1\}$ is partitioned into  $N$-classes of cardinality $ 5 $,
whereas the elements of $\bigcup_{k\in K} ( P^k \setminus \lbrace 1 \rbrace )$ are decomposed into $N$-classes of cardinality 121. Therefore, the set of class sizes of $N$ is $\lbrace 1,  5,  121\rbrace$. Now, let us compute the $G$- class sizes of $N$. As $G$ is a Frobenius group with kernel $K$ and complement {\rm \textbf{N}}$_H(P)$, it follows that $K$ is decomposed exactly into the trivial class and $G$-classes of size $|${\rm \textbf{N}}$_H(P)|= 20$. That is, the $N$-classes contained in $K\setminus \{1\}$ are grouped 4 by 4 to form $G$-classes. And on the other hand, the four $N$-conjugacy classes contained in $\bigcup_{k\in K}P^k \setminus\lbrace 1 \rbrace$ of size 121, are grouped in pairs and become two $G$-conjugacy classes of size $121 \times 2$. Then the set of $G$-class sizes of $N$ is $\lbrace 1, 20,  242\rbrace$ and $\Gamma_G(N)$ is a connected graph.\\

{\bf Example 5.2}  The following example shows that the case in which $N$ is a $p$-group in Theorem E actually occurs. Let $G$ be the group of the library of the small groups of GAP (\cite{gap}) with number Id(324,8) and with the presentation
$$\langle x,y,z|x^3=y^4=z^9=1,[x,y]=1,z^y=z^{-1},z^2=xzxzx=x^{-1}zx^{-1}zx^{-1} \rangle .$$
By using GAP, one can check that $G$ has an abelian normal subgroup $N \cong {\mathbb Z}_3 \times {\mathbb Z}_3$ and the set of $G$-class sizes of $N$ is $\{1,2,3\}$, so $\Gamma_G(N)$ is disconnected.\\

{\bf Open question.} The referee proposed us the following question: whether it is possible to obtain any information on the structure of $G$ from  the graph $\Gamma_G(N)$ or not. We believe that in general  $\Gamma_G(N)$ may provide few information of $G$, although possibly one could get further information on the action of $G$ on $N$. In fact, $G/{\bf C}_G(N)$ is always immerged in Aut$(N)$. For giving an easy example, we consider the case in which $\Gamma_G(N)$ is just one vertex, as it happens with  $G=S_3$ and $N=A_3$. Now, take $N$ any $p$-elementary abelian group of order $p^s$ and let  us consider the action of  $G$= Hol$(N)$ on $N$. As a result of the fact that Aut$(N)$ acts transitively on $N\setminus\{1\}$, it follows that $\Gamma_G(N)$ consists only in one vertex, whilst  Aut$(N)\cong$ GL$(s, p)$ and so $G$ might have a extremely more complex structure.

\bigskip
\noindent {\bf Acknowledgements}

\bigskip

The research of the first and second authors is supported by the Valencian Government,
Proyecto PRO\-METEO/2011/30 and by the Spanish Government, Proyecto
MTM2010-19938-C03-02. The first author is also supported by Universitat Jaume I, grant P11B2012-05.


\begin{thebibliography}{99}

\bibitem{Alf} G. Alfandary.  On graphs related to conjugacy classes of groups. Israel J. Math. 86 (1994), no. 1-3, 211-220.


\bibitem{Akh} Z. Akhlaghi, A. Beltr\'an, M.J. Felipe, M. Khatami, Structure of normal subgroups with three $G$-class sizes. Monatsh. Math. {\bf 167} (2012), 1-12.



\bibitem{BerHerMann} E.A. Bertram, M. Herzog, A. Mann, On a graph related to conjugacy classes of groups. Bull. London Math. Soc. {\bf 22} (6) (1990), 569-575.


\bibitem{Chi} D. Chillag, M. Herzog, A. Mann, On the diameter of a graph related to conjugacy classes of groups. Bull. London Math. Soc. {\bf 25} (1993), 255-262.

\bibitem{Dol} S. Dolfi, Arithmetical conditions on the length of the conjugacy classes of a finite group. J. Algebra {\bf 174} (1995), 753-771.

\bibitem{gap} The GAP Group, GAP - Groups, Algorithms and
Programming, Vers. 4.4.12 (2008). http://www.gap-system.org

\bibitem{Hu} B. Huppert, Character Theory of Finite groups, Vol. 25, De Gruyter Expositions in Mathemathics, Berlin, New York, 1998.






\end{thebibliography}
\end{document}